\documentclass[a4paper,10pt]{article}
\usepackage[english,frenchb]{babel}
\usepackage[T1]{fontenc}
\usepackage[utf8]{inputenc}

\usepackage{amsmath}
\usepackage{array}
\usepackage{amsthm}
\usepackage{amssymb}
\usepackage{setspace}

\usepackage{graphicx}
\usepackage[all]{xy}
\usepackage{geometry}
\input xy

\title{Stabilité homologique pour les groupes d'automorphismes des produits libres}

\author{Ga\"el Collinet\thanks{Institut de Recherche Mathématique avancée (UMR 7501), Université de Strasbourg (collinet@math.unistra.fr).}\,, Aur\'elien Djament\thanks{Laboratoire de mathématiques Jean Leray (UMR 6629), CNRS (aurelien.djament@univ-nantes.fr).}\, et James Griffin\thanks{Department of Pure Mathematics and Mathematical Sciences, University of Cambridge (J.Griffin@dpmms.cam.ac.uk).}}

\newtheorem{th-intro}{Théorème}
\newtheorem{theo}{Théorème}[section]
\newtheorem{pr}[theo]{Proposition}

\newtheorem{lm}[theo]{Lemme}

\theoremstyle{definition}

\newtheorem{nota}[theo]{Notation}

\theoremstyle{remark}
\newtheorem{rem}[theo]{Remarque}
\newtheorem{ex}[theo]{Exemple}

\begin{document}

\maketitle

\begin{abstract}
 On montre dans cet article que, pour tout groupe $G$ indécomposable pour le produit libre $*$ et non isomorphe à $\mathbf{Z}$, l'inclusion canonique ${\rm Aut}(G^{*n})\to {\rm Aut}(G^{* n+1})$ induit un isomorphisme entre les groupes d'homologie $H_i$ pour $n\geq 2i+2$, comme l'avaient conjecturé Hatcher et Wahl. En fait, on montre un peu plus --- en particulier, le résultat vaut pour tout groupe $G$ à condition de remplacer le groupe des automorphismes du produit libre par le sous-groupe des automorphismes symétriques. Nous nous appuyons pour cela sur des constructions et résultats d'acyclicité dus à McCullough-Miller et Chen-Glover-Jensen et sur des propriétés de fonctorialité qui nous permettent d'utiliser des méthodes classiques d'homologie des foncteurs. 

\bigskip

\begin{center}
 {\bf Abstract}
\end{center}

We show in this article that, for any group $G$ indecomposable for the free product $*$ and non-isomorphic to $\mathbf{Z}$, the canonical inclusion ${\rm Aut}(G^{*n})\to {\rm Aut}(G^{* n+1})$ induces an isomorphism between the homology groups $H_i$ for $n\geq 2i+2$, as was conjectured by Hatcher and Wahl. In fact we show a little more --- in particular, the result is true for any group $G$ if we replace the automorphism group of the free product by the subgroup of symmetric automorphisms. For this purpose we use constructions and acyclicity results due to McCullough-Miller and Chen-Glover-Jensen and functoriality properties which allow us to apply classical methods in functor homology.
\end{abstract}

\section*{Introduction}

Soit $G$ un groupe. On note $G^{\ast n}$ le produit libre de $n$ copies de $G$, et
$$\alpha_{n,i}^G: H_i(\mathrm{Aut}( G^{\ast n}))\to H_i(\mathrm{Aut}( G^{\ast n+1})) $$
l'application induite en homologie entière de dimension $i$ par l'inclusion naturelle.

\smallskip
\noindent 
Dans la série d'articles \cite{Hat} -  \cite{HatVog} - \cite{HatVogWah}, A. Hatcher, K. Vogtmann et N. Wahl montrent que $\alpha_{n,i}^{\mathbf{Z}}$ est un isomorphisme pour $n\geq 2i+2$. Les arguments sont d'origine géométrique. En effet, un théorème de F. Laudenbach montre que, pour $V_n=(S^2\times S^1)^{\sharp n}$, l'application $MCG(V_n)\to \mathrm{Aut}(\pi_1(V_n))=\mathrm{Aut}(\mathbf Z^{*n})$ est un épimorphisme dont le noyau est le sous-groupe distingué engendré par les twists de Dehn autour de copies  de $S^2$ dans $V_n$ (la notation $MCG(V)$ désignant le groupe des difféotopies de $V$).

\smallskip
\noindent 
Dans \cite{HatWah}, Hatcher et Wahl montrent des résultats plus généraux sur l'homologie des groupes de difféotopies, et en déduisent que $\alpha_{n,i}^{G}$ est un isomorphisme 
dès que
$n\geq 2i+2$ pour $G=\mathbf Z/n$ avec $n\in\{2,3,4,6\}$ et pour $G=\pi_1(S)$ avec $S$ une surface orientable compacte sans bord, ou $S$ une variété hyperbolique de dimension $3$, de volume fini, n'admettant pas d'isométrie inversant l'orientation.

\smallskip
\noindent Ces auteurs conjecturent alors que $\alpha_{n,i}^{G}$ est un isomorphisme pour $n\geq 2i+2$ quel que soit le groupe $G$.

\medskip Dans cet article, nous démontrons le résultat suivant.

\begin{th-intro}\label{intro1} Soit $G$ un groupe admettant une décomposition $G=G_1\ast \cdots\ast G_p$, chaque $G_i$ étant irréductible pour le produit libre, et aucun des $G_i$ n'étant monogène infini. Alors $\alpha_{n,i}^{G}$ est un isomorphisme pour $n\geq 2i+2$.
\end{th-intro}

En fait, le résultat que nous obtenons est un peu plus général. Pour le décrire et donner une idée de sa démonstration, nous avons besoin d'introduire quelques constructions classiques.

\medskip
\noindent
Soit $\mathrm{G}=(G_e)_{e\in E}$ une famille finie de groupes et $G=\ast_{e\in E} G_e$ leur produit libre. On appelle groupe des automorphismes symétriques de $G$ relativement à $\mathrm{G}$, et on note $\Sigma\mathrm{Aut}(\mathrm{G})$, le sous-groupe de $\mathrm{Aut}(G)$ constitué des automorphismes $\gamma$ tels que pour tout $e\in E$, il existe $e'\in E$ (nécessairement unique) et $g\in G$ tels que l'on ait l'égalité $\gamma(G_e)=gG_{e'}g^{-1}$. \\
L'application $e\mapsto e'$ est alors une permutation de $E$, que l'on note $\pi(\gamma)$, et l'application $\pi:\Sigma\mathrm{Aut}(\mathrm{G})\to \mathfrak{S}(E)$ ainsi définie est un morphisme de groupes dont l'image est notée $\mathfrak S(\mathrm{G})$ et dont le noyau, noté $\mathrm{PAut}(\mathrm{G})$, est appelé groupe des automorphismes purs de $G$ relativement à~$\mathrm{G}$.\\
On se convainc sans peine que pour tout $e\in E$, le noyau de la projection $G\to G_e$ est stable sous l'action de $\mathrm{PAut}(\mathrm{G})$. On en déduit un morphisme $\mathrm{PAut}(\mathrm{G})\to\prod_{e\in E}\mathrm{Aut}(G_e)$ qui est clairement surjectif et scindé. Le noyau de ce morphisme est noté $\mathrm{FR}(\mathrm{G})$ et appelé groupe de Fouxe-Rabinovitch de $G$ relativement à $\mathrm{G}$.\\
Notant $\mathrm{aut}(\mathrm{G})$ le groupe $\prod_{e\in E}\mathrm{Aut}(G_e)\rtimes\mathfrak{S}(\mathrm{G})$, on en déduit une décomposition en produit semi-direct
$$\Sigma\mathrm{Aut}(\mathrm{G})\simeq \mathrm{FR}(\mathrm{G})\rtimes\mathrm{aut}(\mathrm{G})\ \ \ .$$

\medskip
Notre résultat principal prend la forme suivante.

\begin{th-intro}\label{intro2}
 Soit $\mathrm{G}=(G_e)_{e\in E}$ une famille finie de groupes. Soit $M$ un sous-ensemble de $E$, et $\mathrm{H}=(G_e)_{e\in M}$ la sous-famille de $\mathrm{G}$ correspondante. Pour que l'application naturelle $\Sigma\mathrm{Aut}(\mathrm{H})\to \Sigma\mathrm{Aut}(\mathrm{G})$ induise un isomorphisme en homologie entière de dimension $i$ et un épimorphisme en homologie entière de dimension $i+1$, il suffit que 
les deux conditions suivantes soient satisfaites :
\begin{itemize}
 \item l'application naturelle $p : M/\mathfrak{S}(\mathrm{H})\to E/\mathfrak{S}(\mathrm{G})$ est surjective ;
 \item pour toute orbite $U\in M/\mathfrak{S}(\mathrm{H})$ telle que $\mathrm{card}(U)$ et $\mathrm{card}(p(U))$ diffèrent, on a l'inégalité $\mathrm{card}(U)\geq 2i+2$.
\end{itemize}
\end{th-intro}

Esquissons maintenant la démonstration. 

\smallskip
\noindent Le foncteur $\mathrm{G}\mapsto \mathrm{FR}(\mathrm{G})$ défini sur la catégorie discrète dont les objets sont les familles finies de groupes s'étend en un foncteur défini sur la catégorie $\mathbf{Gr}\wr \Theta$ ayant les mêmes objets, et dont les morphismes de source $\mathrm{H}=(H_x)_{x\in X}$ et de but $\mathrm{G}=(G_y)_{y\in Y}$ sont les couples $(f,(\gamma_x)_{x\in s(f)})$ constitués d'une fonction injective $f:X\to Y$ (ici, le terme fonction est utilis\'e avec le sens 
d'{\em application partiellement d\'efinie}) et d'une famille de morphismes $\gamma_x : H_x\to G_{f(x)}$ indexée par le domaine de définition $s(f)$ de $f$.

\smallskip
\noindent Dans la section \ref{s1}, nous définissons, en suivant Eilenberg-MacLane (cf. \cite{EML}, chap.~II), une notion de degré polynomial d'un foncteur $T:\mathbf{Gr}\wr \Theta\to \mathbf{Ab}$, la notation $\mathbf{Ab}$ désignant la catégorie des groupes abéliens. Puis, en s'inspirant de Betley (\cite{Bet}, §\,4), nous montrons un théorème de stabilité pour l'homologie des groupes $\mathrm{Aut}_{\mathbf{Gr}\wr \Theta}(\mathrm{G})$ à coefficients tordus par un foncteur $T$ polynomial. Quelques chasses au diagramme dans des suites spectrales de Lyndon-Hochshild-Serre montrent alors (\ref{clef1}) que pour que le théorème ci-dessus soit vrai, il suffit que le foncteur $\mathrm{G}\mapsto H_i(\mathrm{FR}(\mathrm{G}))$ soit polynomial de degré au plus $2i$ pour tout $i\geq 1$.

\smallskip 
\noindent Cette propriété est vérifiée dans la section \ref{s2}. Nous utilisons une variation de la construction de D. MacCullough et A. Miller \cite{MM}, due à Y. Chen, H. Glover et C. Jensen \cite{CGJ}, d'un espace contractile sur lequel $\mathrm{FR}(\mathrm G)$ agit de manière intéressante. On en déduit (\ref{th-ssh}) une suite spectrale convergeant vers $H_*(\mathrm{FR}(\mathrm G))$ dont la $q$-ème ligne du terme $E^2$ est de la forme 
$H_*(X,A\mapsto H_q(R\circ S))$, où $X$ est un ensemble ordonné fini décrit en termes d'arbres, $S$ est un foncteur de $X$ dans la catégorie des ensembles pointés, et $R$ est un foncteur de la catégorie des ensembles pointés dans les groupes. Cette description nous permet d'utiliser les efficaces techniques d'algèbre homologique sur les petites catégories
(cf. par exemple B. Mitchell \cite{Mit},  T. Pirashvili \cite{Pira-tor}) et J.-L. Loday \cite{Lod}
brièvement rappelées au paragraphe \ref{Ahslpc},  pour montrer que cette suite spectrale s'effondre sur sa première colonne dès le terme $E^2$ et pour identifier précisément cette première colonne (\ref{dec-pira}). La conclusion s'ensuit.

\medskip

Le th\'eor\`eme de Kurosh (cf. par exemple \cite{Serre}, §\,5.5)
implique que si chaque facteur $G_e$ de $\mathrm{G}=(G_e)_{e\in E}$ est irréductible pour le produit libre, et si aucun d'entre eux n'est monogène infini, alors l'inclusion 
$\Sigma\mathrm{Aut}(\mathrm{G})\subset\mathrm{Aut}(\ast_{e\in E} G_e)$ est en fait 
une \'egalit\'e. Le théorème \ref{intro2} implique donc le théorème \ref{intro1}.

\medskip

\noindent
{\em Note.}
Dans un article en préparation \cite{Gri}, le troisième auteur élabore une théorie des {\em produits sur un complexe diagonal}, pour étudier les 
propriétés d'une classe de groupes généralisant les groupes à angles droits d'Artin, qui contient les groupes de Fouxe-Rabinovitch. Deux des résultats 
importants de notre démonstration, les propositions \ref{th-ssh} et \ref{dec-pira}  admettent dans cette théorie un énoncé plus général, 
reposant également sur les constructions de \cite{MM} et \cite{CGJ}, revisitées dans un contexte combinatoire beaucoup plus vaste et adapté au contexte des automorphismes des produits libres. Le phénomène sous-jacent de théorie des groupes est la propriété de 
{\em réductibilité des pics} démontrée par D. Collins et H. Zieschang \cite{ColZie}, qui est l'ingrédient principal aussi bien des démonstrations de 
contractibilité dans \cite{MM} et \cite{CGJ} (contractibilité cruciale dans le présent article) que de la mise en évidence par N. Gilbert \cite{Gil} d'une 
présentation par générateurs et relations de $\mathrm{FR}(\mathrm{G})$ (nécessaire dans \cite{Gri} pour identifier ce groupe à un produit sur un complexe 
diagonal) bien plus simple que celle donnée initialement par
D.I. Fouxe-Rabinovitch \cite{FR}.
Dans \cite{Gri}, le troisième auteur utilise également des constructions introduites par H. Abels et S. Holz dans \cite{AH}. Le résultat final de \cite{Gri} peut permettre de retrouver la stabilité ici abordée ; le but du présent travail consiste à montrer comment on peut se contenter de résultats plus qualitatifs (donc plus faciles à obtenir) que \cite{Gri} pour obtenir la stabilité de l'homologie sans la calculer complètement, à condition de tenir soigneusement compte des fonctorialités. De fait, l'article \cite{Gri} calcule explicitement les foncteurs polynomiaux qui interviennent dans ce travail, dont le caractère qualitatif permet d'espérer des utilisations plus générales pour la stabilité homologique (qu'on ne peut établir qu'exceptionnellement par un calcul complet des groupes d'homologie).

\paragraph*{Remerciements}
 Les deux premiers auteurs témoignent leur gratitude à Nathalie Wahl pour un cours donné à Strasbourg en octobre 2009 qui a attiré leur attention sur le problème de la stabilité homologique pour les automorphismes des produits libres (conjecture d'Hatcher-Wahl). Ils ont également bénéficié du soutien partiel du contrat ANR BLAN08-2\_338236 HGRT {\em Nouveaux liens entre la théorie de l’homotopie et la théorie des groupes et des représentations}.  
 Ils ne soutiennent pas pour autant l'ANR, dont ils revendiquent le transfert des moyens aux laboratoires sous forme de crédits récurrents.

La visite des premier et troisième auteurs au laboratoire de mathématiques Jean Leray (UMR 6629) en février 2011 a été soutenue par la fédération de recherche {\em Mathématiques des Pays de Loire} (FR 2962) à travers le programme Géanpyl.

\section{Effets croisés et stabilité homologique}\label{s1}
\begin{nota} On note
$\mathbf{ens}$ la cat\'egorie des ensembles finis, 
$\mathbf{Gr}$ la cat\'egorie des groupes et 
$\mathbf{Ab}$ celle des groupes ab\'eliens.

\smallskip
La catégorie dont les objets sont les ensembles finis et dont les morphismes sont les fonctions est notée $\Gamma$.  
La sous-catégorie de $\Gamma$ ayant les mêmes objets, mais dans laquelle les morphismes sont les fonctions injectives est notée $\Theta$.

\smallskip Pour un morphisme $f:A\to B$ de l'une de ces catégories, on note $s(f)$ le domaine de définition de $f$. 
\end{nota}

\begin{rem}
 La catégorie $\Gamma$ est équivalente à la catégorie des ensembles
 finis pointés : l'équivalence s'obtient en associant à un ensemble
 pointé le complémentaire de son point de base et  à une application pointée la fonction définie sur le complémentaire de la préimage du point de base du but qu'elle induit.

\end{rem}

\subsection{Effets crois\'es}\label{effcrois}

\begin{rem}
Ce paragraphe s'inspire principalement de \cite{EML}, chap.~II.
\end{rem}

\noindent Dans la cat\'{e}gorie
$\Theta$, on dispose, pour
 tout objet $E$, et tout sous-ensemble $M\subset E$, 
\begin{itemize}
\item[\textbullet] d'un morphisme $i_M^E\in\mathrm{Hom}_{\Theta}(M,
E)$  
repr\'esent\'e par le diagramme $M = M \subset E$,
\item[\textbullet] d'un d'un morphisme
$r_E^M\in\mathrm{Hom}_{\Theta}(E,M) $
 repr\'esent\'e par le diagramme $E\supset M = M$.
\item[\textbullet] et donc d'un idempotent $d_{E,M}=i_M^E\circ r_E^M$ de $E$.
\end{itemize}
On constate que ces idempotents ont les propri\'{e}t\'{e}s suivantes :
\begin{itemize}
\item[P1.] pour tous sous-ensembles $M$ et $N$ de $E$, on a les 
\'{e}galit\'{e}s $d_{E,M}\circ d_{E,N}=d_{E,N}\circ d_{E,M}=d_{E,M\cap N}$,
\item[P2.] pour tout morphisme $f:E\to F$ de $\Gamma$, pour tout
  sous-ensemble $M$ de $E$, on a l'\'{e}galit\'{e} $f\circ
  d_{E,M}=d_{F,f(M)}\circ f$.
\end{itemize} 
 
\medskip

 \noindent Soient $\mathcal{C}$ une cat\'{e}gorie, et $\Lambda$ 
 une sous-cat\'{e}gorie de $\Gamma$ contenant
$\Theta$.
 On note $\mathcal{C}\wr\Lambda$ la cat\'{e}gorie dont 
 \begin{itemize}
 \item[\textbullet] les objets sont les couples $(E, (C_e)_{e\in E})$ 
 compos\'es d'un objet $E$ de $\Lambda$ et d'une famille index\'ee 
 par $E$ d'objets de $\mathcal{C}$,
 \item[\textbullet] les morphismes entre deux objets  $(E, (C_e)_{e\in
 E})$ et $(F, (C_e')_{e\in F})$ sont les couples
 $(f,(\gamma_e)_{e\in s(f)})$ compos\'es d'un morphisme $f:E\to F$ de 
 $\Lambda$ et d'une famille index\'ee par $s(f)$ 
 de morphismes $\gamma_e : C_e\to C_{f(e)}'$ de $\mathcal{C}$.
 \end{itemize}

\noindent Pour un objet $C=(E, (C_e)_{e\in E})$ de $\mathcal{C}\wr\Lambda$ et un 
sous-ensemble $M$ de $E$, on pose $\rho_M(C)=(M, (C_e)_{e\in M})$.
Les morphismes $i_M^E$ et $r_E^M$ de
$\Theta$ se rel\`{e}vent de fa\c{c}on unique en des morphismes 
$i_M^C\in\mathrm{Hom}_{\mathcal C_1\wr\Theta}(\rho_M(C), C)$
et $r_C^M\in\mathrm{Hom}_{\mathcal C_1\wr\Theta}(C,\rho_M(C))$
respectivement, la notation $\mathcal{C}_1$ d\'esignant la cat\'egorie discr\`ete sous-jacente \`a $\mathcal{C}$. 
On note $d_{C,M}= i_M^C\circ r_C^M$ 
l'idempotent associ\'{e} de $\mathrm{End}_{\mathcal C\wr\Lambda}(C)$.
 
\medskip

Soient $\mathcal{A}$ une cat\'{e}gorie ab\'{e}lienne, et 
$T:\mathcal{C}\wr\Lambda\to\mathcal{A}$ un foncteur. 
Pour un objet $C=(C_e)_{e\in E}$ et un sous-ensemble $M$ de
$E$, on pose
$$cr_M(T)(C)=\bigcap_{U\subsetneq M}\mathrm{ker}(T(d_{C,U}))\cap 
\mathrm{im}(T(d_{C,M}))$$
et $cr(T)(C)=cr_E(T)(C)$.

\begin{rem}
 On peut définir de la même façon des effets croisés en supposant que $\Lambda$ contient seulement une sous-catégorie de $\Gamma$ plus petite que $\Theta$ (on a besoin seulement d'inclusions) ; ce raffinement, plus technique à mettre en \oe uvre, ne nous sera pas utile.
\end{rem}

\begin{lm}\label{lmEilMac}
Soient $\Lambda$ une sous-cat\'egorie de $\Gamma$ contenant
  $\Theta$, $\mathcal{C}$ une cat\'egorie, $\mathcal{A}$ une cat\'egorie
  ab\'elienne, et 
$T:\mathcal{C}\wr\Lambda\to\mathcal{A}$ un foncteur tel que
$T(\varnothing)=0$.

\begin{enumerate} 
\item  Soit $f:C\to D$ un morphisme de $\mathcal{C}\wr\Lambda$. 
Alors l'image par $T(f)$ de
$cr_M(T)(C)$ dans $T(D)$ est contenue dans $cr_{f(M)}(T)(D)$ et 
cette image est nulle
si $M$ n'est pas contenu dans $s(f)$.
\item Le morphisme
$T(i_M^C) : cr(T)(\rho_M(C))\to cr_M(T)(C)$ est un isomorphisme.
\item Pour tout objet $C=(C_e)_{e\in E}$  de
$\mathcal{C}\wr\Lambda$, on a l'\'egalit\'e
$$T(C)=\bigoplus_{M\subset E} cr_M(T)(C)\ \ \ .$$
\end{enumerate}
\end{lm}
\begin{proof}
Les deux premiers points sont des cons\'equences imm\'ediates de la
propri\'et\'e P2. D\'emontrons le troisi\`eme point. Pour un ensemble $E$,
notons $\mathcal{P}(E)$ l'ensemble des parties de $E$. 
Soit $C=(C_e)_{e\in E}$ un
objet de $\mathcal{C}\wr\Lambda$, et $T$ un foncteur
$\mathcal{C}\wr\Lambda\to\mathcal{A}$ v\'erifiant l'\'egalit\'e
$T(\varnothing)=0$.
Les endomorphismes $(T(d_{C,M}))_{M\in\mathcal{P}(E)}$ forment d'apr\`es
P1 un
syst\`eme commutatif d'idempotents de $T(C)$. On a donc une \'egalit\'e
$$T(C)= \bigoplus_{P\subset\mathcal{P}(E)}\ \left( \bigcap_{A\in P} 
\mathrm{ker}(T(d_{C,A}))\right) \cap\left( \bigcap_{A\notin P} 
\mathrm{im}(T(d_{C,A})) \right) \ \ \ .$$
Pour $P\subset \mathcal{P}(E)$ fix\'e, posons $M=\cap_{A\in P}A$. Alors
d'apr\`es la propri\'et\'e P1, on a
l'\'egalit\'e 
$$\bigcap_{A\notin P} 
\mathrm{im}(T(d_{C,A}))=\mathrm{im}(T(d_{C,M}))$$
et l'hypoth\`ese $T(\varnothing)=0$ implique que 
$\mathrm{ker}(T(d_{C,A}))\cap \mathrm{im}(T(d_{C,M}))$ est nul d\`es que $A$
n'est pas strictement contenu dans $M$. On en d\'eduit que le sommant
associ\'e \`a  $P\subset\mathcal{P}(E)$ est nul sauf si $P$ est un
  sous-ensemble du type $P=\{ A\in\mathcal{P}(E), A\subsetneq M \} $, d'o\`u
  la conclusion.
\end{proof}

\subsection{Stabilit\'e homologique}\label{stab}

\begin{rem}
Dans ce paragraphe, on s'inspire sans vergogne de travaux de S. Betley. Le théorème~\ref{clef1}
ci-après
généralise en effet le résultat de stabilité homologique à 
coefficients tordus pour les groupes symétriques établi dans le §\,4 de \cite{Bet}. 
\end{rem}
On fixe dans la suite une sous-cat\'egorie $\Lambda$ de $\Gamma$
contenant $\Theta$ et un foncteur $T : \mathcal{C}\wr\Lambda\to\mathbf{Ab}$.

\medskip
\noindent Soit $C=(C_e)_{e\in E}$ un objet de
  $\mathcal{C}\wr\Lambda$. On pose $G_C=\mathrm{Aut}_{\mathcal{C}\wr\Lambda}(C)$. 
On note $V_C$ l'ensemble des \'el\'ements de la forme $\rho_M(C)$, 
$M$ parcourant l'ensemble des parties de $E$.
Un tel \'el\'ement $A=\rho_M(C)$ \'etant donn\'e, on pose $C-A=\rho_{E-M}(C)$. 
On dispose d'un monomorphisme 
$ G_{A}\to G_{C}$ obtenu 
en \'etendant un automorphisme de $A$ par l'identit\'e sur $C-A$, qui fournit un morphisme de $G_A$-modules $T(A)\to \mathrm{res}^{G_{A}}_{G_C}T(C)$ et donc par adjonction un
morphisme de $G_C$-modules $\mathrm{ind}_{G_{A}}^{G_C}T(A)\to T(C)$. On
en d\'eduit un morphisme 
$$\eta_{C,A,i}:H_i(G_{A} ;T(A))
\xrightarrow{_{Shapiro}}H_i(G_C;\mathrm{ind}_{G_{A}}^{G_C}T(A))\to H_i(G_C;T(C)).$$

\begin{pr} On garde les notations  ci-dessus. 
On suppose que :
\begin{enumerate} 
\item $T(\varnothing)=0$,
\item pour $U\in V_C-G_C.V_A$, on a $cr(T)(U)=0$,
\item pour $U\in V_A$ de cardinal
  strictement sup\'erieur \`a $d$ on a $cr(T)(U)=0$,
\item pour $U\in V_A$ de cardinal inf\'erieur ou \'egal \`a $d$, 
l'application $H_q(G_{A-U};\mathbf Z)\to H_q(G_{C-U};\mathbf
  Z)$ est un isomorphisme pour $q\leq i$ et un \'epimorphisme pour
  $q=i+1$.
\end{enumerate}
Alors  
$\eta_{C,A,i} : H_i(G_A;T(A))\to
H_i(G_C;T(C))$
est un isomorphisme et $\eta_{C,A,i+1}$ est un \'epimorphisme.
\end{pr}

La démonstration de cette proposition (ainsi que celle de la proposition
suivante) utilise la remarque suivante.
\begin{rem}\label{remtech} Soient $(E_{p,q}^r,d_{p,q}^r)$ et
  $('\!E_{p,q}^r,'\!d_{p,q}^r) $ suites spectrales du premier quadrant, de type
  homologique (i.e. avec $d^r$ de degré $(-r,r-1)$).
Soit $(\mu_{p,q}^r:E_{p,q}^r\to{}'\!E_{p,q}^r)$ un
  morphisme de suites spectrales. Consid\'erons les propri\'et\'es 
\begin{itemize}
\item[$P_1(r,i)$] le morphisme  $\mu_{p,q}^r$ est un isomorphisme pour tous $p,\ q$ tels que
$p+q\leq i$ et un \'epimorphisme pour tous $p,q$ avec $p+q=i+1$ et
$0<q<i$, 
\item[$P_1'(r,i)$]  le morphisme  $\mu_{p,q}^r$ est un isomorphisme pour tous $p,\ q$ tels que
$p+q\leq i$ et un \'epimorphisme pour tous $p,q$ avec $p+q=i+1$ et
$0\leq q<i$, 
\item[$P_2(r,i)$] la propriété $P_1'(r,i)$ est vérifiée, et le morphisme $\mu_{p,q}^r$ est un \'epimorphisme pour tous $p,\ q$ tels que
$p+q\leq i+1$. 
\end{itemize}

\noindent Le lemme des cinq montre que si $P'_1(r,i)$ est vraie, alors
$P'_1(r+1,i)$ est vraie, et que si $P_2(r,i)$ est vraie, alors
$P_2(r+1,i)$ est vraie. Si de plus on sait que les différentielles
$d_{p,0}^r$ sont nulles pour tout $r$, alors la propriété  $P_1(r,i)$
implique
la propriété $P_1(r+1,i)$.

\noindent Une nouvelle application du lemme des cinq permet d'en déduire que si 
$$\xymatrix{ E^2_{p,q}\ar[d]_{\mu^2_{p,q}}
\ar@{=>}[r] & H_{p+q}\ar[d]^{\mu_{p+q}}\\
'\!E^2_{p,q}\ar@{=>}[r] &
'\!H_{p+q}
}$$
est un morphisme de suites spectrales, alors
\begin{itemize}
\item si la propriété $P_1'(2,i)$ est vérifiée, le morphisme
  $\mu_s:H_s\to\, '\! H_s$ est un isomorphisme pour $s\leq i$,
\item si la propriété $P_2(2,i)$ est vérifiée, le morphisme
  $\mu_s:H_s\to\, '\! H_s$ est un isomorphisme pour $s\leq i$ et un
  épimorphisme pour $s\leq i+1$,
\item Si les différentielles
$d_{p,0}^r$ sont nulles pour tout $r$, et si la propriété $P_1(2,i)$
  est vérifiée, alors $\mu_s:H_s\to\, '\! H_s$ est un isomorphisme pour
  $s\leq i$.
\end{itemize}
\end{rem}

\begin{proof}[Démonstration de la proposition]
On choisit un syst\`{e}me de repr\'esentants $J_C$ du quotient $V_C/G_C$.
Vue dans la cat\'egorie des $G_C$-modules, l'\'egalit\'e du troisi\`eme
point du lemme \ref{lmEilMac} s'\'ecrit   
$$  \bigoplus_{U\in J_C} \mathrm{Ind}_{G_U\times G_{C-U}}^{G_C}(cr(T)(U))\simeq
T(C)\ .$$
Pour $A\in V_C$,  l'application naturelle $V_A
  \to V_C/G_C$ passe au quotient pour donner une injection $V_A/G_A
  \to V_C/G_C$. Posons $J_A=J_C\cap V_A$. L'application $T(i_M^C)$ se d\'ecompose sous la forme
$$ \bigoplus_{U\in J_A } \mathrm{Ind}_{G_U\times G_{A-U}}^{G_A}(cr(T)(U)) 
\to \bigoplus_{U\in
J_C} \mathrm{Res}_{G_C}^{G_A}\mathrm{Ind}_{G_U\times G_{C-U}}^{G_C}(cr(T)(U)).$$ 
Le morphisme $\mathrm{ind}_{G_A}^{G_C}T(A)\to T(C)$ obtenu par
adjonction prend donc la forme 
$$ \bigoplus_{U\in J_A } \mathrm{Ind}_{G_U\times G_{A-U}}^{G_C}(cr(T)(U)) 
\to \bigoplus_{U\in
J_C} \mathrm{Ind}_{G_U\times G_{C-U}}^{G_C}(cr(T)(U)).$$
Le morphisme 
$\eta_{C,A,i} : H_i(G_A;T(A))\to
H_i(G_C;T(C))$
se d\'ecompose ainsi en la somme, sur
$U\in J_A$, des morphismes 
$$ h_{U,i} : H_i(G_U\times G_{A-U};cr(T)(U))\to  
H_i(G_U\times G_{C-U};cr(T)(U))$$
et du morphisme $$r : 0\to \bigoplus_{U\in
J_C-J_A}H_i(G_U\times G_{C-U};cr(T)(U))\ \ \ .$$

Considérons le morphisme de suites spectrales de Lyndon-Hochschild-Serre 
$$\xymatrix{ H_p(G_U;H_q(G_{A-U};cr(T)(U)))\ar[d]
\ar@{=>}[r] & H_{p+q}(G_U\times G_{A-U};cr(T)(U))\ar[d]\\
H_p(G_U;H_q(G_{C-U};cr(T)(U)))\ar@{=>}[r] & H_{p+q}(G_U\times G_{C-U};cr(T)(U))
}\ \ \ .$$
Les actions de $G_U$ sur $H_q(G_{A-U};\mathbf Z)$ et de $G_{C-U}$ sur
$cr(T)(U)$ étant triviales, la
formule des coefficients universels montre que si $H_q(G_{A-U};\mathbf Z)\to H_q(G_{C-U};\mathbf Z)$ est
un isomorphisme pour $q\leq i$ et un \'epimorphisme pour $q=i+1$, la
condition $P_2(2,i)$ de la remarque ci-dessus est vérifiée, et donc que
 $h_{U,i}$ est un isomorphisme et $h_{U,i+1}$ est un
\'epimorphisme. 

Si de plus $cr(T)(U)$ est nul pour $U\in V_A-G(A).V_M$, alors le morphisme
$r$ est un isomorphisme pour tout $i$.
\end{proof}

\medskip
\noindent Un foncteur $T : \mathcal{C}\wr\Lambda\to\mathbf{Ab}$ est dit
polynomial de degr\'e au plus $d$, si $T(\varnothing)=0$ et si $cr(T)(C)$ est nul
 pour tout objet $C=(C_e)_{e\in E}$
avec $\mathrm{card}(E)> d$.

\begin{pr}\label{utile} Soit $\kappa : \mathcal C\wr\Lambda\to
  \mathbf{Gr}$ un foncteur, $C$ un objet de $\mathcal C\wr\Lambda$. Si 
\begin{itemize}
\item pour $q>0$, le foncteur $H_q\circ\kappa$ est polynomial de
degr\'e au plus $d(q)$,
\item l'inclusion $V_A/G_A\subset V_C/G_C$ est une \'egalit\'e,
\item pour tout $q\in [1,i-1]$, pour tout $U\in V_A$ non vide de cardinal
inf\'erieur ou \'egal \`a $d(q)$, le morphisme 
$$ \sigma_{\ell,U} : H_\ell(G_{A-U};\mathbf Z)\to H_\ell(G_{C-U};\mathbf Z)$$
est un isomorphisme pour $\ell\leq i-q$ et un \'epimorphisme pour
$\ell = i-q+1$,
\item $\sigma_{\ell,\varnothing}$ est un isomorphisme pour $\ell\leq i$,
\end{itemize}
alors  le morphisme 
$$\mu_\ell\ :\ H_\ell(\kappa(A)\rtimes G_A;\mathbf Z)
\to H_\ell(\kappa(C)\rtimes G_C;\mathbf Z)$$ est un isomorphisme
  pour $\ell\leq i$.
Si de plus 
\begin{itemize} 
\item $\sigma_{i+1,\varnothing}$ est un \'epimorphisme,
\item pour tout $U\in V_A$ non vide de cardinal
inf\'erieur ou \'egal \`a $d(i)$, le morphisme 
$ \sigma_{1,U} $ est un \'epimorphisme, 
\end{itemize}
alors $\mu_{i+1}$ est un \'epimorphisme.
\end{pr}

\begin{proof}
Notons $E^r_{p,q}(C)$ le terme g\'en\'eral de la suite spectrale 
de Lyndon-Hochschild-Serre associ\'ee \`a l'extension
$G_C\ltimes \kappa(C)$, et pour $A\in V_C$, consid\'erons le morphisme
de suites spectrales
$$\xymatrix{ E^2_{p,q}(A)=H_p(G_A;H_q(\kappa(A)))\ar[d]_{\mu^2_{p,q}}
\ar@{=>}[r] & H_{p+q}(G_A\ltimes \kappa(A);\mathbf Z)\ar[d]^{\mu_{p+q}}\\
E^2_{p,q}(C)=H_p(G_C;H_q(\kappa(C)))\ar@{=>}[r] &
H_{p+q}(G_C\ltimes \kappa(C);\mathbf Z)
}$$
Les extensions \'etant scind\'ees, les morphismes $d^r_{p,0}$ sont
nuls. 
Il suffit donc, d'après la remarque~\ref{remtech}, pour pouvoir conclure
que $\mu_\ell$ est un isomorphisme pour $\ell\leq i$ (resp. un
\'epimorphisme pour $\ell\leq i+1$), de v\'erifier
que la propri\'et\'e $P_1(2,i)$ (resp. $P_2(2,i)$) est v\'erifi\'ee. C'est
exactement ce que l'on peut d\'eduire des hypoth\`eses gr\^ace \`a la proposition pr\'ec\'edente. 
\end{proof}

\medskip

Soit $C=(C_x)_{x\in X}$ un objet de $\mathcal{C}\wr\Theta$. La relation 
$\sim$ d\'efinie sur $X$  par $x\sim x'$ si $C_x$ et $C_x'$ sont isomorphes 
(dans $\mathcal{C}$) est une relation d'\'equivalence. 
On note $C(x)$ la classe d'\'equivalence de $x$.
Choisissons un syst\`eme de repr\'esentants $S$ de $X/\sim$. 
On a un isomorphisme 
 $$\mathrm{Aut}_{\mathcal{C}\wr\Theta}(C)\simeq 
\prod_{x\in S}\left( \left(\prod_{x'\in
C(x)}\mathrm{Aut}_{\mathcal{C}}(C_{x'})\right)
\rtimes\mathfrak{S}(C(x))\right)\ .$$
Notons $G_C$ ce groupe d'automorphismes.

\begin{theo}\label{clef1}
Soit $\kappa : \mathcal{C}\wr\Theta\to\mathbf{Gr}$ un foncteur tel que 
pour tout $q>0$, le foncteur $C\mapsto H_q(\kappa(C);\mathbf Z)$ 
soit polynomial de degr\'e au plus $2q$.
Soient $C=(C_x)_{x\in X}$ un objet de $\mathcal{C}\wr\Theta$, $M$ un 
sous-ensemble de $X$, et $A=(C_x)_{x\in M}$. 
Pour que le morphisme naturel 
$$\eta_{C,A,\ell}: H_\ell(\kappa(A)\rtimes G_A;\mathbf Z)\to 
H_\ell(\kappa(C)\rtimes G_C;\mathbf Z)$$
soit un isomorphisme pour $\ell\leq i$ et un \'epimorphisme pour $\ell\leq i+1$, 
il suffit que
\begin{itemize}
\item l'inclusion naturelle $M/\sim\to X/\sim$ soit une \'egalit\'e,
\item pour tout $x\in M$, l'inclusion $A(x)\to C(x)$ soit 
une \'egalit\'e ou bien que le cardinal de $A(x)$
soit sup\'erieur 
ou \'egal \`a $2i+2$.
\end{itemize}
\end{theo}

\begin{proof} Elle se fait en deux \'etapes :

\noindent A] On consid\`ere d'abord la situation suivante : 
\begin{itemize}
\item $\mathcal{C}$ est la cat\'egorie discr\`ete $\{1,\dots,m\}$,
\item $K_1,\dots,K_m$ sont des groupes,
\item $\kappa : \mathcal{C}\wr\Lambda\to\mathbf{Gr}$ est le foncteur qui \`a $(n_i)_{i\in I}$ associe
$\prod_{i\in I}K_{n_i}$.
\end{itemize}
Le th\'eor\`eme de Künneth montre que le foncteur $H_q\circ\kappa$ est de degr\'e au plus $q$.

\noindent Les travaux de Nakaoka \cite{Nak} sur l'homologie des groupes sym\'etriques nous
apprennent que si $Y$ est un ensemble fini, si $X$ est un sous-ensemble de $Y$, alors l'application 
$H_p(\mathfrak{S}(X))\to H_p(\mathfrak{S}(Y))$ est un monomorphisme pour tout $p$, et est  un
isomorphisme si $\mathrm{card}(X)\geq 2p$.

\noindent On d\'eduit donc de la proposition \ref{utile} que si $C$ est un objet de 
$\mathcal{C}\wr\Lambda$, si $A$ est un \'el\'ement de $V_C$, alors l'application 
$$ H_\ell\left( \prod_{n=1}^m K_n\wr\mathfrak{S}(A(n))\right)\to 
H_\ell\left( \prod_{n=1}^m K_n\wr\mathfrak{S}(C(n))\right)$$
est un isomorphisme pour $\ell\leq i$ (resp. un \'epimorphisme pour $\ell \leq i+1$) d\`es que l'on a l'\'egalit\'e $V_A/G_A=V_C/G_C$,
et que pour tout $n\in [1,m]$, ou bien l'inclusion $A(n)\subset C(n)$ est une \'egalit\'e, 
ou bien on a $\mathrm{card}(A(n))\geq 2i+1$ (resp. $\mathrm{card}(A(n))\geq 2i+2$),

\noindent B] Le th\'eor\`eme s'ensuit alors directement de A] et de la proposition \ref{utile} appliqu\'ee \`a la situation suivante :
\begin{itemize}
\item  $\mathcal{C}$ est une cat\'egorie quelconque,
\item  $\kappa : \mathcal{C}\wr\Lambda\to\mathbf{Gr}$ est un foncteur tel que pour tout $q$, 
le foncteur $H_q\circ \kappa$ soit de degr\'e au plus $2q$.
\end{itemize}
\end{proof}

\subsection{Algèbre homologique sur les petites catégories}\label{Ahslpc}

On utilisera l'écriture $a\in J$ pour spécifier que $a$ est un objet d'une petite catégorie~$J$. 

\medskip
\noindent
Soient $k$ un anneau commutatif (dans la suite de ce texte, on prendra toujours $k=\mathbf Z$), et $J$ une petite catégorie.
On note
\begin{itemize} 
\item $k[J]$ l'anneau dont le groupe abélien sous-jacent est 
$$\bigoplus_{a,b\in J}k[\mathrm{Hom}_J(a,b)]$$
 et dont la multiplication est définie par 
$$\gamma.\gamma'=
\left\{\begin{array}{ccl} \gamma\circ\gamma' &\mathrm{si} & \mathrm{source}(\gamma)=\mathrm{but}(\gamma')\\ 0 & \mathrm{sinon}&
\end{array}\right.\ \ \ ,$$
\item $J-\mathrm{Mod}$ la catégorie des foncteurs $J\to k-\mathrm{Mod}$,
\item $P_a^J:=k[\mathrm{Hom}_{J}(a,-)]$ pour $a\in J$,
\item $\underline J=\bigoplus_{a\in J}P_a^J$,
\item $\underline{\underline J}: J^{op}\times J\to k-\mathrm{Mod}$  le bifoncteur $k[\mathrm{Hom}_{J}(-,-)]$.
\item $\underline k$ le foncteur identiquement égal à $k$.
\end{itemize} 

\medskip
\noindent
Un résultat bien connu de Mitchell \cite{Mit} §\,7
affirme d'une part que le foncteur évident $\varepsilon : J-\mathrm{Mod}\to k[J]-\mathrm{Mod}$ est pleinement fidèle et d'autre part que
si $J$ n'a qu'un nombre fini d'objets
alors $\varepsilon$ est une équivalence de catégories. Dans ce travail, on pourrait se ramener, au prix de quelques spécialisations, à cette situation, 
mais nous préférons faire des rappels d'ordre plus général.

\medskip
\noindent
La catégorie $J-\mathrm{Mod}$ hérite de la catégorie $k[J]-\mathrm{Mod}$, via $\varepsilon$,
les
propriétés suivantes :
\begin{itemize}
\item c'est une catégorie abélienne,
\item elle est munie d'un produit tensoriel $-\otimes_{J}-: J^{op}-\mathrm{Mod}\times J-\mathrm{Mod}\to k$,
\item elle a assez de projectifs : l'égalité 
$\varepsilon(\underline J)=k[J]$ montre que les foncteurs $P_a^J$ sont projectifs (ce qui s'ensuit aussi du lemme de Yoneda), 
et que tout foncteur de $J-\mathrm{Mod}$  admet une résolution projective (et admet même une résolution par des foncteurs qui sont 
somme directe de copies de $\underline J$).
\end{itemize} 

\medskip
\noindent
Soient $F\in J-\mathrm{Mod}$ et $G\in J^{op}-\mathrm{Mod}$ deux foncteurs.
On note les isomorphismes canoniques
$$\underline{k}\otimes_J F\simeq\mathrm{colim}_J F\ \ \ ,\ \ \ \mathrm{Hom}_{J-\mathrm{Mod}}(P_a^J,F)\simeq F(a)\ \ \ ,
\ \ \ P_a^{J^{op}}\otimes_J F\simeq F(a)\ \ \ ,\ \ \ G \otimes_J P_a^{J}\simeq G(a)\ \ \ .$$

\medskip
\noindent
Soient $C_*$ et $D_*$ deux complexes de
chaînes
$\mathbf N$-gradués de $J-\mathrm{Mod}$ et $J^{op}-\mathrm{Mod}$ respectivement. 
On a deux suites spectrales $^IE^1_{p,q}=H_p(H_q(C_*)\otimes_J D_*)$ et $^{II}E^1_{p,q}=H_p(C_*\otimes_J H_q(D_*))$ qui convergent vers
$H_{p+q}(Tot(C_*,D_*))$,
la notation $Tot(C_*,D_*)$ désignant le complexe total du double complexe $C_*\otimes_J D_*$. On en déduit que le calcul
des foncteurs dérivés de $-\otimes_{J}-$, notés $\mathrm{Tor}^J_*$ peut se faire indifféremment en résolvant la première ou la seconde variable, et
que l'on a de manière générale une suite spectrale 
\begin{eqnarray}\label{suitspecpira} E^2_{p,q}=\mathrm{Tor}^J_p(F,H_q(C_*))& \Rightarrow &H_{p+q}(F\otimes_J C_*)\ \ \ .\end{eqnarray}
 
\bigskip
Soient $J$ et $K$ deux petites catégories et $F:J\to K$ un foncteur.  
Par précomposition de $\underline{\underline K}$ avec $F$, on en déduit un bifoncteur 
$F^*(\underline{\underline K})\in J^{op}\times K-\mathrm{Mod}$. Le foncteur 
$F_!(-)=F^*(\underline{\underline K})\otimes_J -: J-\mathrm{Mod}\to K-\mathrm{Mod}$ 
est adjoint à gauche du foncteur $F^*:K-\mathrm{Mod}\to J-\mathrm{Mod}$, et l'on a un isomorphisme naturel
\begin{eqnarray}\label{adjkan} F^*(M)\otimes_J N\simeq M\otimes_K F_!(N)\ \ \ .\end{eqnarray}

\medskip
\noindent
Soient $N\in J-Mod$ et $M\in K-\mathrm{Mod}$. 
De (\ref{suitspecpira}) et (\ref{adjkan}), on déduit une suite spectrale 
\begin{eqnarray}\label{suitspecgroth} E^2_{p,q}=\mathrm{Tor}^K_p(M,\mathbf{L}_qF_!(N))\Rightarrow \mathrm{Tor}^J_{p+q}(F^*(M),N) \ \ \ ,\end{eqnarray}
la notation $\mathbf{L}_qF_!(N)$ désignant le $q$-ème foncteur dérivé à gauche de $F_!$ évalué en $N$. Ces foncteurs dérivés sont naturellement isomorphes 
à $ \mathrm{Tor}^J_{q}(F^*(\underline{\underline K}),N) : a\mapsto 
\mathrm{Tor}^J_{q}(F^*(P_a^{K^{op}}),N)$. 

On peut encore écrire $\mathbf{L}_qF_!(N)(a)\simeq H_q(J^a_F,\pi^*N)$ naturellement où $J^a_F$ est la catégorie dont les objets sont les couples $(b,\varphi)$ où $b\in J$
 et $\varphi\in {\rm Hom}_K(F(b),a)$, les morphismes $(b,\varphi)\to (b',\varphi')$ étant les morphismes $u : b\to b'$ de $J$ tels que $\varphi=\varphi'\circ F(u)$, et $\pi : J^a_F\to J\quad (b,\varphi)\mapsto b$ le foncteur canonique. (En effet, on dispose d'une bijection canonique
$${\rm Hom}_J(b,\pi(x))\simeq\bigsqcup_{\varphi\in {\rm Hom}_J(F(b),a)}{\rm Hom}_{J^a_F}((b,\varphi),x),$$
d'où l'on déduit $\pi^*P^J_b\simeq\bigoplus_{\varphi\in {\rm Hom}_J(F(b),a)} P^{J^a_F}_{(b,\varphi)}$ puis 
$${\rm Tor}^{J^a_F}_*(G,\pi^*F)\simeq{\rm Tor}^J\Big(b\mapsto\bigoplus_{\varphi\in {\rm Hom}_J(F(b),a)}G(b,\varphi),F\Big)\quad.)$$

En particulier, prenant $N=\underline{k}$, notre suite spectrale (\ref{suitspecgroth}) peut s'écrire :
\begin{eqnarray}\label{groth} E^2_{p,q}=\mathrm{Tor}^J_p(M,a\mapsto H_q(J^a_F))& \Rightarrow &\mathrm{Tor}^J_{p+q}(F^*(M),\underline k)=H_{p+q}(J,F^*(M))\ \ \ .\end{eqnarray}

\begin{ex}
 Soit $i:\Omega\to \Gamma$ l'inclusion de la sous-catégorie de $\Gamma$ ayant les mêmes objets et dont les morphismes sont les surjections partout définies. On se permet de remplacer implicitement ces deux catégories par des squelettes pour en faire de petites catégories. Pour $T\in \Omega-\mathrm{Mod}$ et $a\in \Gamma$, on a
$$i_!(T)(a)= i^*(P_{a}^{\Gamma^{op}})\otimes_\Omega T=\left(\bigoplus_{E\in \Omega} k[\mathrm{Hom}_\Gamma(i(E),a)]\otimes_k T(E)\right)/S\ \ \ ,$$
le sous-module $S$ étant engendré par les éléments $\varphi\otimes x - id_M\otimes T(\varphi)(x)$ avec $\varphi\in\mathrm{Hom}_\Omega(M',M)$ et $x\in T(M')$.
On en déduit un isomorphisme
$$i_!(T)(E)=\bigoplus_{M\subset E}T(M)\ \ \ .$$
Le premier point du lemme \ref{lmEilMac} montre que les effets croisés définissent un foncteur $cr:\Gamma-\mathrm{Mod} \to
  \Omega-\mathrm{Mod} $, et le second et le troisième point de ce même lemme montrent que $i_!$ est un inverse de $cr$. 
Cette équivalence de Morita montre en particulier que l'on a des isomorphismes naturels
\begin{eqnarray}\label{piraeq} \mathrm{Tor}_*^\Gamma(M,N)\simeq\mathrm{Tor}_*^\Omega(cr(M),cr(N))\ \ \ .\end{eqnarray}
Cette équivalence est due à Pirashvili (cf. \cite{PDK}, ou \cite{PH}).
Pour des résultats d'équivalence de Morita~\guillemotleft~à la Dold-Kan~\guillemotright~extrêmement généraux (mais n'utilisant pas explicitement d'effets croisés), on pourra aussi se reporter à l'article \cite{Slo} de S{\l}ominska.
\end{ex}

\section{Homologie des groupes de Fouxe-Rabinovitch et des automorphismes symétriques d'un produit libre}\label{s2}
Soient $\mathrm{G}=(G_e)_{e\in E}$ une famille finie de groupes et $G$ le produit libre de cette famille.
Le groupe de Whitehead $\mathrm{Wh}(\mathrm{G})$ de $\mathrm{G}$ est la préimage par l'application naturelle $\mathrm{PAut}(\mathrm{G})\to\prod_{e\in E}\mathrm{Aut}(G_e)$ du sous-groupe
 $\prod_{e\in E}\mathrm{Inn}(G_e)$,
la notation $\mathrm{Inn}(H)$ désignant le groupe des automorphismes intérieurs d'un groupe $H$. Il est engendré par les automorphismes du type 
$$\alpha^{g_i}_{i,j} : \gamma_e\in G_e \mapsto \left\{\begin{array}{ccc} \gamma_e & \mathrm{si} & e\neq j\\ g_i\gamma_eg_i^{-1} & \mathrm{si} & e = j\end{array}\right.  \ \ \mathrm{avec}\ \  g_i\in G_i \ \ \ .$$
On note les relations 
\vskip -1.5cm
\begin{spacing}{2}
  \begin{eqnarray}\label{relFR}
 \begin{array}{lclcc}
  \alpha^{g_i}_{i,j}\circ \alpha^{g_i'}_{i,j} &=&\alpha^{g'_ig_i}_{i,j}&\mathrm{si}&\mathrm{card} \{i,j\}=2\\
\alpha^{g_i}_{i,j}\circ \alpha^{g_i'}_{i,j'} &=&\alpha^{g_i}_{i,j'}\circ \alpha^{g_i'}_{i,j}&\mathrm{si}&\mathrm{card} \{i,j,j'\}=3\\
 \alpha^{g_i}_{i,j}\circ (\alpha^{g_i'}_{i',j}\circ \alpha^{g_i'}_{i',i})&=&(\alpha^{g_i'}_{i',j}\circ \alpha^{g_i'}_{i',i})\circ\alpha^{g_i}_{i,j} &\mathrm{si}&\mathrm{card}\{i,i',j'\}=3.
 \end{array}
\end{eqnarray}
 \end{spacing}
\vskip -0.5cm

Ce groupe de Whitehead contient en particulier le sous-groupe $\mathrm{Inn}(G)$, qui est distingué. Le quotient $\mathrm{Wh}(\mathrm{G})/\mathrm{Inn}(G)$ est noté $\mathrm{OWh}(\mathrm{G})$.
 
\bigskip
Dans \cite{MM}, D. McCullough et A. Miller construisent un ensemble ordonné contractile $\mathrm{MM}(\mathrm{G})$ de dimension $\mathrm{card}(E)-2$ sur lequel agit de façon intéressante le groupe $\mathrm{OWh}(\mathrm{G})$ : (a)~l'action est fidèle, (b) le quotient est compact et (c) les stabilisateurs des drapeaux sont bien compris.

\smallskip
\noindent Posons $E_+=E\sqcup\{\star\}$ et complétons G en $\mathrm{G}_+=(G_e)_{e\in E_+}$ en  choisissant un groupe non trivial $G_\star$.
Dans \cite{CGJ}, Y. Chen, H. Glover et C. Jensen constatent que l'injection $\mathrm{Wh}(\mathrm{G})\to \mathrm{OWh}(\mathrm{G_+})$ permet 
de construire un $\mathrm{Wh}(\mathrm{G})$-sous-ensemble ordonné du complexe $\mathrm{MM}(\mathrm{G}_+)$ ayant les propriétés (a), (b) et (c), 
et démontrent la contractibilité de ce nouvel ensemble ordonné (qui ne dépend pas du choix de $G_\star$) pour $G_{\star}=\mathbf Z/2$.

\noindent Soit $\mathbf{Gr/Po}$
la catégorie des actions de groupes sur des ensembles ordonnés,
autrement dit la 
catégorie 
\begin{itemize}
 \item dont les objets sont les 
triplet $A=(G_A,X_A,\mu_A)$ consitués d'un groupe $G_A$, 
d'un ensemble ordonné $X_A$, et d'un morphisme 
$\mu_A : G_A\to \mathrm{Aut}(X_A)$, 
\item dont les morphismes de source un objet $A$ et de but 
un objet $B$ sont les couples $(\gamma,\varphi)$ constitués 
d'un morphisme de groupes $\gamma : G_A\to G_B$ et d'un morphisme 
d'ensembles ordonnés $\varphi:X_A\to X_B$ tels que l'on ait les égalités 
$\varphi(\mu_A(g).x)=\mu_B(\gamma(g)).\varphi(x)$ quels que soient $g\in G_A$ et $x\in X_A$. 
\end{itemize}

\smallskip
\noindent L'objet de cette section 
est de montrer comment utiliser la construction de \cite{CGJ} pour obtenir un foncteur
$\mathcal{C}:\mathbf{Gr}\wr\Theta\to\mathbf{Gr/Po}$
et d'en tirer les conclusions qui s'imposent.

\subsection{Le complexe de Chen-Glover-Jensen}
Nous proposons ici un bref apperçu de la construction. Pour plus de détails, le lecteur peut se référer à  \cite{MM} et \cite{CGJ}.

\smallskip
\noindent Soit $J_E$ l'ensemble des arbres bipartites
\'{e}tiquet\'{e}s par $E$, plantés en $\star$, ordonné par la relation de pliage (voir la figure 1). 
Dans un tel arbre, il y a trois types de sommets : des sommets étiquetés, des sommets muets, et la racine $\star$. 
Un sommet étiqueté n'a pour voisins que des sommets muets, un sommet muet n'a pas de voisin muet, et n'est jamais une feuille,
la racine a un unique voisin muet. Le plus petit élément de $J_E$, qui ne contient qu'un sommet muet, 
et dont tous les sommets étiquetés sont des feuilles, est noté $N_E$.

\begin{figure}[!h]\label{pliage}
\centering
\includegraphics[height=4cm]{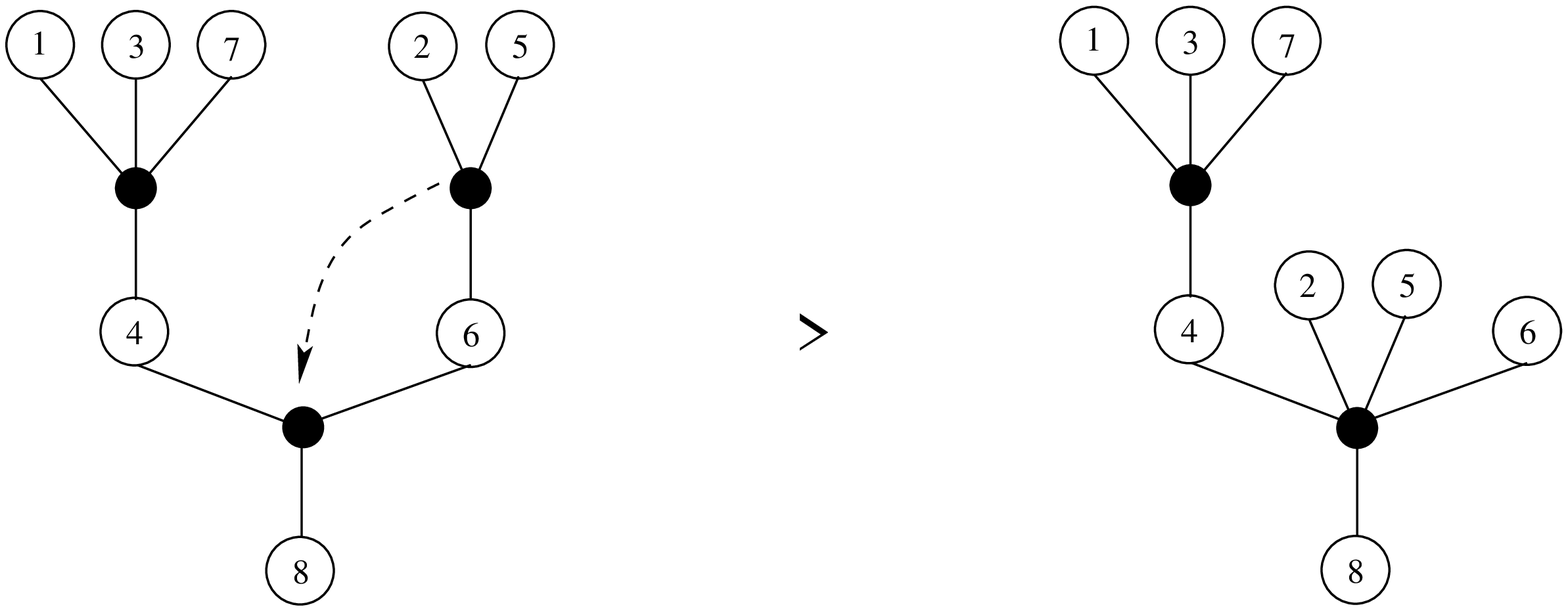}
\caption{Un pliage d'arbres bipartites plantés étiquetés par $\{1,\dots,7\}$}
\end{figure}

\smallskip
\noindent  Soient $\mathrm{G}=(G_e)_{e\in E}$ une famille finie de groupes et $A$ un élément de $J_{E}$. Un automorphisme de Whitehead de $G$
relatif \`{a} la base $\mathrm{G}$, de facteur op\'{e}rant
$\ell\in E$  est support\'{e} par $A$ si, pour toute paire de sommets étiquetés $(i,j)$ contenus dans 
une m\^{e}me composante connexe de $A-{\ell}$, on a l'\'{e}galit\'{e} $g_{\ell,i}=g_{\ell,j}$, et si $g_{\ell,i}$ est trivial lorsque $i$ est 
dans la composante connexe de $\star$.

\smallskip
\noindent Soit $\mathcal B$ l'ensemble des familles $\rho(\mathrm{G})=(\rho(G_1),\dots,\rho(G_n))$ obtenues en appliquant à $\mathrm{G}$ un élément de $\mathrm{Wh}(\mathrm{G})$.
Soit $\sim$ la relation d'\'{e}quivalence sur le  produit $\mathcal{B}\times
J_E$ engendrée par
 $(\mathrm{B},A)\sim(\rho\mathrm{B},A)$ 
pour tout automorphisme
de Whitehead $\rho$ de $G$ relatif \`{a} $\mathrm{B}$ support\'{e} par $A$. 

\smallskip 
\noindent Notation : si $\mathrm{B}$ est un élément de $\mathcal{B}$ et $A$ un élément de $J_E$, on note $(\underline{\mathrm{B}},A)$ la classe d'équivalence de $(\mathrm{B},A)$.

\noindent 
La relation d'ordre~$\leq$ d\'{e}finie plus haut s'\'{e}tend au quotient $(\mathcal{B}\times
\mathcal{A})/\sim$ :
$$(\underline{\mathrm{B}},A)\leq(\underline{\mathrm{B}}',A')\ \ \ \Leftrightarrow 
\ \ \ (\mathrm{B},A')\sim  (\mathrm{B}',A')\ \ \ \mathrm{et}\ \ \ A\leq A'$$
et en fait un ensemble ordonn\'{e} que nous notons $X_{\mathrm G}$, 
muni d'une action (par automorphismes d'ensembles ordonnés) du groupe
$\mathrm{Wh}(\mathrm{G})$ (définie sur les représentants par $\alpha.(\mathrm{B},A)=(\alpha.\mathrm{B},A)$).

\medskip

\noindent {\bf Th\'eor\`eme de Chen-Glover-Jensen.} {\em L'ensemble ordonné $X_{\mathrm{G}}$ est contractile.}

\begin{rem} Dans \cite{CGJ}, on suppose d\`es le d\'{e}but
que tous les groupes $G_i$ sont finis, mais cette hypoth\`ese
n'intervient pas dans la d\'{e}monstration du premier th\'{e}or\`eme reformul\'{e} ci-dessus.
\end{rem}

\subsection{Action du groupe de Fouxe-Rabinovitch sur $X_{\mathrm{G}}$}

On garde les notations du paragraphe pr\'{e}c\'{e}dent. On note $E_{\mathrm{G}}$ l'étoile du sommet $(\underline{\mathrm{G}},N_E)$.

\begin{lm} L'action du groupe 
$\mathrm{FR}(\mathrm{G})$ sur
$X_\mathrm{G}$ admet pour domaine fondamental le sous-ensemble ordonné $E_{\mathrm{G}}$.
\end{lm}

\begin{proof} Il suffit de démontrer que l'action est 
transitive sur les sommets minimaux. 
Disons que deux sommets minimaux sont
voisins si l'intersection de leurs \'{e}toiles est non vide. Soit
$({\mathrm{H}},N)$ un  sommet minimal de $X_\mathrm{G}$. 
Il existe une
suite de sommets minimaux $({\mathrm{H}}_i,N)_{i=1}^k$ tels que l'on
ait
$$({\mathrm{H}}_0,N)=(\mathrm{G} ,N)\ \ \ ,\ \ \ ({\mathrm{H}}_k,N)=(\mathrm{H} ,N)$$
et tels que pour $i\in [k-1]$, le sommet $(\mathrm{H}_i ,N)$ soit
voisin du sommet $(\mathrm{H}_{i+1} ,N)$.\\
Il suffit donc de d\'{e}montrer que deux sommets minimaux voisins 
sont dans la même orbite.
Si $(\mathrm{H},N)$ et
$(\mathrm{H}',N)$ sont voisins, il existe un arbre $A$ tel que les
sommets $(\underline{\mathrm{H}},A)$ et $(\underline{\mathrm{H}}',A)$
soient \'{e}gaux, autrement dit il existe une suite
$(\rho_1,\dots,\rho_s)$ d'automorphismes de Whitehead support\'{e}s par
$A$, relatifs \`{a} la base $\mathrm{H}$, tels que l'on ait l'\'{e}galit\'{e} 
$\rho_s\cdots\rho_1.\mathcal{G}=\mathcal{H}$. Il suffit donc, pour
que l'affirmation soit vraie, qu'elle le soit pour $s=1$, 
ce qui est \'{e}vident.
\end{proof}
                       
\begin{nota} On note 
$\alpha_G:\mathrm{FR}(\mathrm{G})\to \mathrm{Aut}(X_G)$ 
l'action du groupe de Fouxe-Rabinovitch de 
$\mathrm{G}$ sur $X_{\mathrm{G}}$.
\end{nota}

\subsection{Fonctorialité}

La construction précédente permet d'associer à tout objet $\mathrm{G}=(G_e)_{e\in E}$ de $\mathbf{Gr}\wr\Theta$ un objet 
$\mathcal C(\mathrm{G})=(\mathrm{FR}(\mathrm{G}),X_{\mathrm{G}},\alpha_{\mathrm{G}})$ de $\mathbf{Gr/Po}$. Montrons comment étendre $\mathcal C$ en un foncteur. 

\medskip 
\noindent On commence par étendre $E\mapsto J_E$ en un foncteur $J:\Theta\to \mathbf{Po}$. Soit $f:E\stackrel{u}{\hookleftarrow} M\stackrel{v}{\hookrightarrow} E'$ un morphisme de $\Theta$. On note $J_f:J_E\to J_{E'}$ le morphisme d'ensembles ordonnés associant à un arbre planté $A$ étiqueté par une fonction $s:E\to \mathrm{sommets}(A)$ l'arbre planté $A'$ étiqueté par la fonction $s':E'\to \mathrm{sommets}(A')$ obtenu par la recette suivante :
\begin{itemize} 
 \item pour chaque sommet $S$ de $A$ étiqueté par un élément de $E-u(M)$, on replie sur une même arête toutes les arêtes issues de $S$,
\item on coupe toutes les feuilles étiquetées par un élément de $E-u(M)$, 
\item on étiquette  le résultat $A''$ par la fonction $s':v(M)\to \mathrm{sommets}(A'')$ définie par $s'(v(m))=s(u(m))$,
\item on ajoute, pour chaque élément $e$ de $E'-v(M)$, une feuille, étiquetée par $s'(e)$, voisine de l'unique voisin muet de $\star$.
\end{itemize}

\medskip
\noindent Soient maintenant $E$ et $E'$ deux ensembles finis, $\mathrm{G}=(G_e)_{e\in E}$ et $\mathrm{G}'=(G_e')_{e\in E'}$ deux objets de $\mathbf{Gr}\wr\Theta$
 et $f:\mathrm{G}\to\mathrm{G}'$ un morphisme de $\mathbf{Gr}\wr\Theta$. On note $\gamma:\mathrm{FR}(\mathrm{G})\to\mathrm{FR}(\mathrm{G}')$ le morphisme qui s'en déduit. 
Soit $\mathrm{B}=\alpha(\mathrm{G})$ un élément de $\mathcal{B}_\mathrm{G}$ (avec $\alpha\in \mathrm{FR}(\mathrm{G})$). 
On pose $\varphi_0(\mathrm{B})=\gamma(\alpha).\mathrm{G}'$.
On définit $\varphi_f:X_{\mathrm{G}}\to X_{\mathrm{G}'}$ par la formule $\varphi_f((\underline{\mathrm{B}},A))=(\underline{\varphi_0(\mathrm{B})},J_f(A))$.

\medskip
\noindent Le foncteur $\mathcal{C}$ auquel il est fait allusion dans l'introduction de cette section est obtenu en posant $\mathcal{C}(f)=(\gamma_f,\varphi_f)$.

\subsection{Stabilisateurs}

Soit $A$ un objet de $J_E$ (que l'on considère indifféremment comme un
ensemble ordonné ou une catégorie). Pour $e\in E$, notons $S_e^A$ le sommet étiqueté par
$e$, et $\mathcal{P}_e^A$ l'ensemble des composantes connexes de $A-S_e^A$ distinctes de celle de la racine $\star$.

\medskip
\noindent Soit $A\leq A'$ un morphisme de $J_E$. La partition
$\mathcal{P}_e^{A'}$ raffine la partition $\mathcal{P}_e^{A}$. Soit
$f_e^{A,A'}$ la fonction associant à un élément de l'ensemble
$\mathcal{P}_e^{A'}$ (une partie) l'élément le contenant de
$\mathcal{P}_e^{A}$ (cette fonction n'est donc définie que sur les parties de $\mathcal{P}_e^{A'}$ qui ne sont pas contenues dans la composante connexe de la racine de $A$).

\medskip
\noindent Soit $F_E : J_E\to (\Gamma^E)^{op}$ le foncteur
associant
\begin{itemize}
\item à un arbre $A$ de $J_E$ l'objet $(\mathcal{P}_e^A)_{e\in E}$
  (chaque partition $\mathcal{P}_e^A$ étant naturellement pointée par
  la partie contenant la racine),
\item à un morphisme $A\leq A'$ le morphisme $(f_e^{A,A'})_{e\in E}$.
\end{itemize}

\medskip
\noindent Soient maintenant $\mathrm{G}=(G_e)_{e\in E}$ un objet de
$\mathbf{Gr}\wr\Theta$. On note $D_{\mathrm{G}}$ le foncteur
$({\Gamma^E})^{op}\to\mathbf{Gr}$ associant à $(P_e)_{e\in E}$ le groupe
$\prod_{e\in E}G_e^{P_e-\{\star\}}$.

\medskip
\noindent Soit $(\underline{\mathrm{G}},A)$ un sommet de
$X_{\mathrm{G}}$ contenu dans l'étoile de $(\mathrm{G},A)$. 
Le stabilisateur $\mathcal{S}(\underline{\mathrm{G}},A)$ de $(\underline{\mathrm{G}},A)$ sous l'action de $\mathrm{FR}(\mathrm{G})$ est 
le sous-groupe formé des automorphismes de Fouxe-Rabinovitch supportés par $A$ relativement à la famille $\mathrm{B}$. 
Les relations \ref{relFR} montrent que l'application
$$\begin{array}{ccc}
D_G\circ F_E (A) &\to& \mathcal{S}(\underline{\mathrm{B}},A)\\
((g_{P})_{P\in \mathcal{P}_{e}^A-\{\star\}})_{e\in E}&\mapsto&  \underset{e\in E}{\prod}\ \underset{P\in \mathcal{P}_e^A}{\prod}\ \underset{p\in P}{\prod} \alpha_{e,p}^{g_P}
\end{array}$$ est un isomorphisme.

\medskip
\noindent En particulier, le stabilisateur d'un drapeau de $X_{\mathrm{G}}$ est aussi celui de son plus petit sommet, et pour 
$(\underline{\mathrm{B}},A)\leq (\underline{\mathrm{B}},A')$, l'injection $\mathcal{S}(\underline{\mathrm{B}},A)\to \mathcal{S}(\underline{\mathrm{B}},A')$
applique un générateur $\prod_{p\in P} \alpha_{e,p}^{b_P}$ sur le
produit $\prod_{P'\in\mathcal{P}_e^{A'},\ P'\subset P} \prod_{p\in
  P'}\alpha_{e,p}^{b_{P}}$, autrement dit, est $D_\mathrm{G}\circ F_E
(A\leq A')$.

\medskip Puisque $X_{\mathrm G}$ est contractile, la construction de
  Borel fournit finalement le résultat suivant.

\medskip 
\begin{pr}\label{th-ssh}  Les foncteurs $\mathrm{FR}:\mathbf{Gr}\wr\Theta\to \mathbf{Gr}$ {\em et}
 $\mathrm{G}=(G_e)_{e\in
  E}\mapsto\mathrm{colim}_{J_E}{D_\mathrm{G}\circ F_E}$  sont isomorphes. 
 De plus, il existe une suite spectrale naturelle en
 $\mathrm{G}=(G_e)_{e\in E}\in \mathbf{Gr}\wr\Theta$ 
$$E^2_{p,q}=H_p(J_E;\mathcal{H}_q(D_\mathrm{G}\circ F_E;\mathbf{Z}))\Rightarrow
H_{p+q}(\mathrm{FR}(\mathrm{G});\mathbf Z)$$
(où $\mathcal{H}_q(D_\mathrm{G}\circ F_E;\mathbf{Z})$ désigne le foncteur $A\mapsto H_q(D_\mathrm{G}\circ F_E(A);\mathbf{Z})$).
\end{pr}

\subsection{Effondrement de la suite spectrale}

Soit $\mathcal{F}_E$ le sous-ensemble de $J_E$ constitué des arbres $A$
dont tous les sommets muets, sauf peut-être le sommet voisin de la
racine, sont bivalents.

Le résultat suivant est essentiellement une reformulation du théorème~3.5 de \cite{Gri}.

\begin{pr}\label{dec-pira} Soit $ :
  (\Gamma^E)^{op}\to\mathbf{Ab}$ un foncteur. Il existe un isomorphisme naturel
 $$H_0(J_E;T\circ F_E)\simeq\bigoplus_{A\in\mathcal{F}_E} cr_{F_E(A)}(T)$$
 tandis que $H_n(J_E;T\circ F_E)$ est nul pour $n>0$.
\end{pr}

\begin{proof}
 On examine la suite spectrale  \ref{groth}
$$E^2_{p,q}={\rm Tor}^{\Gamma^E}_p(T,X\mapsto H_q(I_X))\Rightarrow H_{p+q}(J_E;T\circ F_E)$$

la notation $I_X$ désignant, pour tout ensemble fini $X$, l'ensemble ordonné des couples $(A,\varphi)$ composés d'un
arbre $A$ de $J_E$ et d'une fonction $\varphi=(\varphi_e)_{e\in E}\in \Gamma^E(F_E(A),X)$, la
relation d'ordre étant définie par
$$(A,\varphi)\leq (A',\varphi')
\text{ si }A\leq A'\text{ et }\varphi'=\varphi\circ F_E(A\leq A').$$

Soit $(A,\varphi)$ un élément de $I_X$, et $e$ un sommet étiqueté de $A$, à distance $k$ de la racine. Les composantes connexes de $A-e$ sont en bijection 
avec l'ensemble $M_e$ des arêtes $(e,m)$ où $m$ parcourt l'ensemble des voisins de $e$, qui sont tous à distance $k+1$ de la racine, 
sauf un, disons $m_e$, qui est à distance $k-1$. 
Soit $\rho_e(A,\varphi)$ l'unique élément de $I_X$ inférieur ou égal à $(A,\varphi)$ dont l'arbre est obtenu en repliant sur $(e,m_e)$ toutes les arêtes $(e,m)$ de $M_e$ 
sur lesquelles $\varphi_e$ n'est pas définie. 
L'application $\rho_e$ est un endomorphisme idempotent de l'ensemble ordonné $I_X$ vérifiant l'inégalité $\rho_e(x)\leq x$ pour tout $x$, 
et fournit donc une rétraction de $I_X$ sur son image. 
De plus, pour tous $e$ et $e'$, les rétractions $\rho_e$ et $\rho_e'$ commutent. Soit $\rho^{\downarrow}$ la composée de toutes ces rétractions. 
L'image $I'_X$ de $\rho^{\downarrow}$ est constituée des couples 
$(A,\varphi)$  dans lesquels $\varphi$ est une {\em application}.
Soit $ (A,\varphi)$ un tel élément. Pour un sommet étiqueté $e$ à distance $k$ de la racine,  notons $\rho^e(A,\varphi)$ l'unique élément de $I'_X$ supérieur 
ou égal à 
$(A,\varphi)$ dans lequel tous les sommets muets voisins de $e$ à distance $k+1$ sont bivalents.  
L'application $\rho^e$ est un endomorphisme idempotent de l'ensemble ordonné $I'_X$ vérifiant l'inégalité $\rho^e(x)\geq x$ pour tout $x$, 
et fournit donc une rétraction de $I'_X$ sur son image. De plus, pour tous $e$ et $e'$, les rétractions $\rho_e$ et $\rho_e'$ commutent. 
Soit $\rho^{\uparrow}$ la composée de toutes ces rétractions. L'image de $\rho^{\uparrow}$ étant discrète,
les composantes connexes de
$I_X$ sont contractiles, et on a  un isomorphisme 
$H_0(I_X)\simeq\underset{A\in\mathcal{F}_E}{\bigoplus}
  R_A(X)$,
la notation $R_A$ désignant le foncteur $X\mapsto R_A(X)=\mathbf{Z}[\mathbf{ens}^E(F_E(A),X)]$.
En particulier, notre suite
spectrale se r\'eduit \`a un isomorphisme naturel
$$H_n(J_E;T\circ F_E)\simeq\bigoplus_{A\in\mathcal{F}_E}{\rm Tor}^{\Gamma^E}_n( T,R_A)\ \ .$$
L'isomorphisme \ref{piraeq} montre que, pour tout objet $Y$ de $\Gamma^E$, le foncteur
$\mathbf{Z}[\mathbf{ens}^E(Y,-)]: (\Gamma^E)^{op}\to\mathbf{Ab}$ est
projectif et que son produit tensoriel
au-dessus de $\Gamma^E$ avec $T$ s'identifie naturellement à $cr_Y(T)$.
\end{proof}

\begin{rem}\label{rq-fonct}

 On peut préciser complètement la fonctorialité en $E$ de l'isomorphisme donné par la proposition~\ref{dec-pira}. 
Le cas d'un morphisme général de $\Theta$ nécessite quelques précautions ; toutefois il est facile de voir qu'une injection (partout définie) 
d'ensembles finis $E\hookrightarrow E'$ induit l'inclusion canonique entre les facteurs correspondants via la proposition --- 
noter que l'inclusion ensembliste $\mathcal{F}_E\hookrightarrow\mathcal{F}_{E'}$ induite par notre injection permet d'identifier 
$F_E(A)$ et $F_{E'}(A)$ pour $A\in\mathcal{F}_E$.
\end{rem}

\subsection{Conclusion}
\begin{lm}\label{polypoly} Soit $T :(\Gamma^E)^{op}\to\mathbf{Ab}$ un foncteur
  polynomial de degré au plus $q$. Le foncteur $\Theta\to\mathbf{Ab}$
  associant à un ensemble $E$ le groupe abélien
 $\bigoplus_{A\in\mathcal{F}_E} cr_{F_E(A)}(T)$
est polynomial de degré au plus $2q$.
\end{lm}

\begin{proof}
Soit $A$ un élément de $\mathcal{F}_E$. Le cardinal de $F_E(A)$ est le nombre $M'$ de sommets muets de $A$ moins 1. 
Soient $F_0$ le nombre de feuilles voisines de l'unique voisin muet de la racine, 
$F_1$ le nombre des autres sommets étiquetés voisins de ce même sommet, $F_2$ le nombre de tous les sommets restant.
On a l'égalité $F_0+F_1+M'=\mathrm{card}(E)$ et l'inégalité $F_1\leq M'$, d'où l'on déduit l'inégalité $F_0\geq \mathrm{card}(E)-2M'$. On en déduit que pour 
$\mathrm{card}(E)\geq 2q+1$, l'application 
$$\bigoplus_{M\subsetneq E}\bigoplus_{A\in\mathcal{F}_M} cr_{F_M(A)}(T)\to \bigoplus_{A\in\mathcal{F}_E} cr_{F_E(A)}(T)$$ est surjective.
\end{proof}

Le théorème \ref{th-ssh}, la proposition \ref{dec-pira}, et le lemme
\ref{polypoly} mis ensemble fournissent enfin le résultat suivant.
\begin{pr}\label{degpol} Le foncteur $\mathrm{G}\mapsto
  H_q(\mathrm{FR}(\mathrm{G});\mathbf Z)$ 
est polynomial de degré au plus $2q$.
\end{pr} 
Combinée au théorème~\ref{clef1}, la proposition~\ref{degpol} fournit le théorème~\ref{intro2}, et donc le théorème~\ref{intro1} qui en constitue un cas particulier, annoncés en introduction.

\bibliographystyle{smfalpha}
\bibliography{biblio-prodlib}

\end{document}